\documentclass[11pt]{article}

\usepackage{amsfonts,amsmath,graphicx}
\graphicspath{{.}}
   \makeatletter
   \let\accent@spacefactor\relax
   \makeatother
   \usepackage[french]{babel}
\usepackage{amsfonts}
\usepackage{amsmath}   %multiline etc. Subsumes amstext,amsbsy,amsopn
\usepackage{amsthm}    %provide proof env
\usepackage{amscd}     %commutative diagrams
\usepackage{euscript}

\def\P{{\Bbb P}}

\def\I{{\Bbb I}}

\newtheorem{defi}{D\'{e}finition}[section]
\newtheorem{pro}[defi]{Proposition}

\newtheorem{lem}[defi]{Lemme}
\newtheorem{theo}[defi]{Th\'{e}or\`{e}me}

\newtheorem{rem}[defi]{Remarque}

\title{
{\bf Fibr\'es logarithmiques sur le  plan projectif}}
\author{Jean Vall\`es}
\date{16 juillet 2005}

\begin{document}

\maketitle
{\small \textbf{R\'esum\'e :} Nous d\'ecrivons le sch\'ema des droites de saut des fibr\'es logarithmiques sur le plan projectif (thm \ref{principal} de ce  texte). Connu, depuis l'article \cite{DK} de Dolgachev et Kapranov pour les fibr\'es de premi\`ere classe de Chern paire, ce r\'esultat  est nouveau lorsque la premi\`ere classe de Chern est impaire.} 

\smallskip

{\small \textbf{Abstract :} We describe   the scheme of jumping lines of logarithmic vector bundles on the projective plane (thm \ref{principal} of this text).  
This result is already proved by Dolgachev and Kapranov in \cite{DK} when the first Chern class is even, it is new when the first Chern class is odd.}

\smallskip

{ \footnotesize ``Ainsi les clart\'es, en s'accumulant, font figure d'enigmes, \`a la mani\`ere d'un verre trop \'epais qui cesse d'\^etre transparent.''(Simone Weil)}
%%%%%%%%%%%%%%%%%%%%%%%%%%%%%%%%%%%%%%%%%%%%%%%%%%%%%%%%%%%%
\section{Introduction}
%%%%%%%%%%%%%%%%%%%%%%%%%%%%%%%%%%%%%%%%%%%%%%%%%%%%%%%%%
Dolgachev et Kapranov ont \'etudi\'e les fibr\'es logarithmiques, introduits par Deligne, dans le cas particulier o\`u la base est un espace projectif complexe $\P^{n}$ et le diviseur associ\'e une r\'eunion d'hyperplans
en position lin\'eaire g\'en\'erale (voir les articles \cite{De} et \cite{DK}) . Ils montrent  notamment que la donn\'ee du diviseur d\'etermine le fibr\'e logarithmique associ\'e (th\'eor\`eme 7.2 dit de "Torelli", voir aussi \cite{Va} pour une preuve diff\'erente). 

\smallskip

Sur $\P^{2}$ les fibr\'es logarithmiques sont stables de rang deux. Selon que leur premi\`ere classe de Chern est paire ou impaire, leur sch\'ema de droites de saut est, en g\'en\'eral, une courbe pour le cas pair, un sch\'ema de longueur fini pour le cas impair. Dans les deux cas le degr\'e de la courbe ou la longueur du sch\'ema fini sont connus et donn\'es en fonctions des classes de Chern. 

\smallskip

Puisque, d'apr\`es le th\'eor\`eme dit de "Torelli",  le fibr\'e est d\'etermin\'e par la donn\'ee du diviseur (ici une r\'eunion de droites), ce diviseur doit aussi d\'eterminer le sch\'ema des droites de saut du fibr\'e associ\'e. Le lien est clairement \'etabli dans \cite{DK} (voir
thm 5.2)
pour les fibr\'es de premi\`ere classe de Chern paire via les courbes mono\"{\i}dales
mais il n'est pas fait pour les fibr\'es de premi\`ere classe de Chern impaire. Plus pr\'ecis\'ement, bien que les droites du diviseur 
soient compt\'ees avec multiplicit\'e,
elles ne suffisent  pas \`a remplir le sch\'ema de droites de saut. Il y a d'autres droites dans le sch\'ema. 

\smallskip

On se propose dans cet article, \`a partir d'un lien simple entre un fibr\'e 
logarithmique sur le plan projectif et le faisceau d'id\'eaux du groupe de points
associ\'e dans le plan projectif dual, de dire qui sont ces droites suppl\'ementaires (i.e nous donnons  la description du sch\'ema des droites de saut des fibr\'es logarithmiques de premi\`ere classe de Chern impaire, th\'eor\`eme \ref{principal}). Ainsi, de la m\^eme fa\c{c}on que  dans le cas pair 
une unique courbe appel\'ee courbe mono\"idale \'etait associ\'ee aux droites de d\'epart (le diviseur du fibr\'e logarithmique),  dans le cas impair un autre ensemble fini de droites sera associ\'e aux droites de d\'epart.

\smallskip

L'\'etude  des fibr\'es vectoriels stables via leur sch\'ema des droites de saut est souvent ingrate dans le sens o\`u les justifications g\'eom\'etriques (pourquoi cette droite saute t-elle et pas une autre?) n'apparaissent pas toujours. De ce point de vue l'int\'er\^et des fibr\'es logarithmiques est que les justifications g\'eom\'etriques sont faciles \`a mettre en \'evidence. C'est ce que nous faisons dans la derni\`ere partie en explicitant la   g\'eom\'etrie sous-jacente du r\'esultat principal de cet article. 
%%%%%%%%%%%%%
\subsection{Notations}
%%%%%%%%%%%%%%%%%%%%%%%
On note $\I \subset \P^2 \times \P^{2 {\vee}}$ la
vari\'et\'e d'incidence
 points-droites de $\P^2$,  $p$ et $q$ les projections  $\P^2
\stackrel{ p}\leftarrow \I \stackrel{ q}\rightarrow
\P^{2\vee}$.

\smallskip

On note $M_{\P^2}(c_{1},c_{2})$ l'espace de modules des fibr\'es vectoriels semi-stables sur $\P^{2}$ de classes de Chern $c_{1}$ et $c_{2}$. En g\'en\'eral on pr\'ef\`erera travailler avec les fibr\'es normalis\'es, c'est \`a dire dont la premi\`ere classe de Chern vaut $0$ ou $-1$. On notera si n\'ec\'essaire ${\cal E}^{\rm norm}$ le fibr\'e normalis\'e correspondant \`a ${\cal E}$. 

\smallskip

\'Etant donn\'e un fibr\'e stable ${\cal E}$ de rang deux on note  $S({\cal E})$ son
sch\'ema des droites de saut et $S_{\textrm{ens}}({\cal E})$ l'ensemble de ses droites de saut. 
%%%%%%%%%%%%%%%%%%%%%%%%%%%%%%%%%%%%%%%%%%%%%%%%%%%%%%%%%%%%%%%%
\section{Lien entre sections globales et droites de saut}
%%%%%%%%%%%%%%%%%%%%%%%%%%%%%%%%%%%%%%%%%%%%%%%%%%%%%%%%%%%%
Soit $F$ un fibr\'e stable de rang deux sur $\P^2$ tel que 
$F\in M_{\P^2}(0,2n-\epsilon)$ avec $\epsilon = 0, 1$. On rappelle qu'une droite 
$l\in \P^{2\vee}$ est une droite de saut pour $F$ ssi 
$F_l=O_l(a)\oplus
O_l(b)$ avec $\mid a-b\mid >1$. On suppose qu'il n'existe pas de droite $l$ telle que 
$H^0F_l(-2)\neq 0$ (cas g\'en\'eral). L'id\'ee principale est alors de comparer le fibr\'e 
$F$ et le fibr\'e $q_{*}p^{*}F$. Comme pour $l$ g\'en\'erale 
 $F_l=O_l\oplus O_l$, on a une application 
$$ F \mapsto (q_{*}p^{*}F)^{\rm norm} \in M_{\P^{2\vee}}(0,n(n-1)) \quad {\rm si }\,\,
\epsilon=0$$
$$ F \mapsto (q_{*}p^{*}F)^{\rm norm} \in M_{\P^{2\vee}}(-1,(n-1)^2) \quad {\rm si }\,\,
\epsilon=1$$
De plus le fibr\'e $q_{*}p^{*}F$ est un fibr\'e de Steiner (voir par exemple \cite{DK} pour la terminologie),
plus pr\'ecis\'ement il s'inscrit dans une suite exacte 
$$ 0 \rightarrow q_{*}p^{*}F \longrightarrow H^1(F(-1))\otimes O_{\P^{2\vee}}(-1)
\longrightarrow H^1(F)\otimes O_{\P^{2\vee}} \rightarrow 0 $$
 On a alors une premi\`ere proposition (ensembliste) liant un
fibr\'e et son image directe.
%%%%%%%%%%%%
\begin{pro}\label{lien}
$$S_{\textrm{ens}}(q_{*}p^{*}F)=\{x\in \P^2 \mid h^0(F\otimes {\frak m}_{x}^{a}(a))\neq 0 \,\, {\rm pour }
\,\,a<[\frac {2n-\epsilon}{2}] 
\}$$
\end{pro}
%%%%%%%%%%%%%%%%%%
{\bf Preuve de la proposition \ref{lien}.} On note  $\P^{2 }(x)$ l'\'eclatement de $\P^{2 }$ le long
d'un point $x\in \P^{2 }$ et   $p_{x}$ et $ q_{x}$  les morphismes
apparaissant dans le diagramme d'incidence
$$
\begin{CD}
 \P^{2 }(x) @>q_x>>    x^{ \vee}\\
 @Vp_xVV \\
\P^{2 }
\end{CD}
$$
On v\'erifie facilement que l'absence de bisauteuse implique
$q_{x*}p_{x}^{*}F=(q_{*}p^{*}F)_{\mid x^{\vee}}$.

\smallskip

Soit $x\in \P^2$, on suppose que la droite correspondante $x^{\vee}$ est une droite
de saut de $q_{*}p^{*}F$. On a alors la d\'ecomposition  
$q_{*}p^{*}F_{\mid x^{\vee}}=O_{x^{\vee}}(-a)\oplus
O_{x^{\vee}}(-b)$ avec $a+b=2n-\epsilon$ et $a<[\frac {2n-\epsilon}{2}]$. Comme 
$q_{x*}p_{x}^{*}F=(q_{*}p^{*}F)_{\mid x^{\vee}}$, cette d\'ecomposition
fournit une section non nulle
$h^0(p_{x}^{*}F\otimes q_{x}^{*}O_{x^{\vee}}(a))\neq 0$ ce qui \'equivaut
\`a $h^0(F\otimes {\frak m}_{x}^{a}(a))\neq 0 $.

\smallskip

R\'eciproquement, soit $s \in H^0(F\otimes {\frak m}_{x}^{a}(a))$ une section non nulle et 
$a<[\frac {2n-\epsilon}{2}]$.  Elle induit une suite  exacte
$$ 0 \rightarrow O_{\P^2} \longrightarrow F(a) \longrightarrow {\cal J}_{\Gamma}(2a) 
\rightarrow 0 $$ o\`u $\Gamma$ est un sh\'ema de longueur finie v\'erifiant   ${\cal
J}_{\Gamma}\subset {\frak m}_{x}^{a}$.

\smallskip

 La section 
$p_{x}^{*}O_{\P^2} \longrightarrow p_{x}^{*}F(a)$ s'annule en codimension 1 le long de
$a$-fois le diviseur exceptionnel $p^{-1}(x)$. Apr\`es simplification on obtient une
section $q_{x}^{*}O_{x^{\vee}}(-a) \longrightarrow p_{x}^{*}F$ ce qui prouve que 
$O_{x^{\vee}}(-a) \hookrightarrow q_{*}p^{*}F_{\mid x^{\vee}}$. $\Box$
%%%%%%%%%%
%%%%%%%%%%%%%%%%%%%%%%%%%%%%%%%%
\section{Fibr\'es logarithmiques du plan}
%%%%%%%%%%%%%%%%%%%%%%%%%%%%%%%%%%%%%%%%%%%%%%%%%%%%%%%%%%%
Soit $Z$ un groupe de points de longueur $2n+1-\epsilon$ en position lin\'eaire
g\'en\'erale de
$\P^2$. C'est une cons\'equence du th\'eor\`eme 7.2 de \cite{DK} que le fibr\'e logarithmique $E(Z)$ sur $\P^{2 \vee}$
associ\'e
\`a $Z$ est d\'efini et enti\`erement d\'etermin\'e par les donn\'ees (a1) et (a2): 

\smallskip

a1) $E(Z)$ est un fibr\'e de Steiner appartennant \`a $M_{\P^{2\vee}}(0,n(n-1))$ si
$\epsilon = 0$, \\ 
appartennant \`a $M_{\P^{2\vee}}(-1,(n-1)^2)$ si
$\epsilon = 1$.

\smallskip

a2) Pour $x\in Z$ on a   
$E(Z)_{\mid x^{\vee}} = O_{x^{\vee}}\oplus O_{x^{\vee}}(2n-\epsilon)$
%%%%%%%%%%
\begin{theo} 
\label{principal}
a) $E(Z)^{\vee}(-1)=q_{*}p^{*}{\cal J}_{Z}(1) $\\
b)  Lorsque $\epsilon =0$,  $S(E(Z))= \overline{\{x\in \P^{2}, h^0({\cal J}_{Z}\otimes {\frak
m}_x^{n-1}(n))\neq 0 \}}$.\\
c) Lorsque $\epsilon =1$, $S(E(Z))$ est support\'e par le lieu g\'eom\'etrique $$ Z \cup \overline{\{x\in \P^{2}, h^0({\cal J}_{Z}\otimes {\frak
m}_x^{n-2}(n-1))\neq 0 \}}$$ et si $Z$ est  g\'en\'eral parmi les groupes de points du plan de longueur $2n$ alors
$$S(E(Z))= Z^{n-1}  \sqcup\{x\in \P^{2}, h^0({\cal J}_{Z}\otimes {\frak
m}_x^{n-2}(n-1))\neq 0 \}$$ o\`u $Z^{n-1}$ d\'esigne le $(n-2)$-i\`eme voisinage infinit\'esimal de $Z$ et $\sqcup$ l'union disjointe.
\end{theo}
%%%%%%%%%%
{\bf Remarque 1.} La courbe $\overline{\{x\in \P^{2}, h^0({\cal J}_{Z}\otimes {\frak
m}_x^{n-1}(n))\neq 0 \}}$ est not\'ee $C(Z)$ par Dolgachev et Kapranov et est appel\'ee
complexe mono\"{\i}dal associ\'e \`a $Z$. Le point $b)$ de ce th\'eor\`eme est d\'ej\`a
prouv\'e, avec une plus grande g\'en\'eralit\'e, par ces deux auteurs (thm 5.2,
\cite{DK}). \\
{\bf Remarque 2.} Lorsque $Z$ est en position g\'en\'erale un simple calcul de dimension du syst\`eme lin\'eaire permet de v\'erifier que $S(E(Z))$ est fini lorsque $\epsilon =1$ (c'est pourquoi il n'est plus n\'ec\'essaire de prendre la cl\^oture de Zariski de l'ensemble). \\
{\bf Remarque 3.} Par contre lorsque $\epsilon =0$, l'irr\'eductibilit\'e de la courbe 
mono\"{\i}dale associ\'e \`a $Z$, m\^eme lorsque $Z$ est en position g\'en\'erale n'est pas 
\'etablie (sauf pour $n\le 3$).

\smallskip

{\bf Preuve du th\'eor\`eme \ref{principal}.} On consid\`ere une extension g\'en\'erale 
$$ 0 \rightarrow O_{\P^2} \longrightarrow F(1) \longrightarrow {\cal J}_Z(2) \rightarrow
0$$ Comme $Z$ est en position lin\'eaire g\'en\'erale il n'existe pas de tris\'ecante \`a
$Z$, ce qui \'equivaut au fait qu'il n'existe pas de droite $l$ telle que $h^0F_l(-2)\neq
0$. On en d\'eduit alors que le fibr\'e $q_{*}p^{*}F=q_{*}p^{*}{\cal J}_{Z}(1)$ est un
fibr\'e de Steiner avec les bonnes classes de Chern. Si $x\in Z$ il existe une section non
nulle 
$ s\in H^0(F\otimes {\frak m}_x(1))$, ce qui implique (comme dans la preuve de la
proposition pr\'ec\'edente) $ q_{*}p^{*}F_{\mid x^{\vee}}= 
O_{x^{\vee}}(-1)\oplus O_{x^{\vee}}(2n+1-\epsilon)$, ce qui prouve $a)$.

\medskip

Maintenant si $x\notin Z$ on a, pour tout entier $a$, une suite exacte 
$$ 0 \rightarrow {\frak m}_x^{a}(a-1) \longrightarrow F\otimes {\frak m}_x^{a}(a) 
\longrightarrow {\cal J}_Z\otimes {\frak m}_x^{a}(a+1)
\rightarrow 0$$
De plus  $H^0F\otimes {\frak m}_x^{n-1}(n-1)=H^0 {\cal J}_Z\otimes {\frak m}_x^{n-1}(n)$
et $H^0F\otimes {\frak m}_x^{n-2}(n-2)=H^0 {\cal J}_Z\otimes {\frak m}_x^{n-2}(n-1)$. 
Compte tenu de la proposition \ref{lien} ceci prouve $b)$ et $c)$. 

\smallskip

Pour prouver le point $d)$, nous calculons maintenant la longueur 
du sous-sch\'ema de saut support\'e par le diviseur $Z$. 
%%%%%%%%%%%%%%%%%%%%%%
\begin{lem}\label{singodd}
Soit $l$ une droite de saut d'un fibr\'e $E\in M_{\P^2}(-1,n)$. Si
$h^0(E_l(-k))=1 $ avec $k>0$ alors $S(E)$
contient le  $(k-1)$-i\`eme voisinage infinit\'esimal du point  $l\in
\P^{2\vee}$.
\end{lem}
%%%%%%%%%%%%%%%%%%%%%
{\bf Preuve du lemme \ref{singodd}.}  L'hypoth\`ese fournit un
homomorphisme surjectif de faisceaux
$$E \longrightarrow O_l(-k-1) $$
Notons $F(-1)\in M_{\P^2}(-2,n-k)$ le fibr\'e noyau de cet homorphisme. 
En appliquant le foncteur $q_{*}p^{*} $ on
trouve
$$ R^1q_{*}p^{*}F(-1)\longrightarrow R^1q_{*}p^{*}E \longrightarrow R^1q_{*}p^{*}O_l(-k-1) \rightarrow 0$$
Tout d'abord, on rappelle que $R^1q_{*}p^{*}O_l(-k-1)=O_{\P^{2\vee}}/{\frak
m}_l^{k}$ (il suffit, pour le prouver de  prendre 
dans $\P^2$ une r\'esolution de $O_l(-k-1)$ et d'utiliser 
 $R^1q_{*}p^{*}O_{\P^2}(-m)=S^{m-2}(\Omega (1))$). La surjection 
$$R^1q_{*}p^{*}E \longrightarrow O_{\P^{2\vee}}/{\frak m}_l^{k} \rightarrow 0$$
implique alors (voir \cite{GP} prop. B2 et B3 par exemple) que $S(E)$ contient le 
$(k-1)$-i\`eme voisinage infinit\'esimal du point  $l\in
\P^{2\vee}$. Enfin, si
$l\notin \textrm{supp}(R^1q_{*}p^{*}F(-1))$ i.e. si $l$ n'est pas une
droite de saut de $F$ alors le
sous-sch\'ema de  $S(E)$ support\'e par $l$ est exactement le
$(k-1)$-i\`eme voisinage infinit\'esimal de $l$
et il est de longueur $\binom{k+1}{2}$.~$\Box$
%%%%%%%%%%%%%%%%%%%%

\smallskip 

Revenons au cas des fibr\'es logarithmiques. Supposons donn\'e $Z=\lbrace l_{1},\cdots, l_{2n}\rbrace$ notre groupe de points 
 en position g\'en\'erale.  On a une suite exacte 
canonique (par \cite{DK} prop. 2.9)
$$\begin{CD}
0 @>>> E(\lbrace l_{1},\cdots, l_{2n-1}\rbrace) @>>> E(\lbrace l_{1},\cdots, l_{2n}\rbrace) @>>>  O_{l_{2n}}
 @>>> 0 
\end{CD}$$
Comme la premi\`ere classe de Chern de $E(\lbrace l_{1},\cdots, l_{2n-1}\rbrace)$ est paire son sch\'ema de droites de saut est une courbe. Si l'on   choisit $ l_{2n}$
hors du support de cette courbe, le sous-sch\'ema de 
$S(E(\lbrace l_{1},\cdots, l_{2n}\rbrace))$ support\'e par $ l_{2n}$ est exactement le 
$(n-2)$-i\`eme voisinage infinit\'esimal de $ l_{2n}$. On en d\'eduit  que 
le sous-sch\'ema support\'e par $Z$ est exactement $Z^{n-1} $. De plus comme la 
la longueur du  sch\'ema de droites de saut de $E(Z)$ est \'egale \`a $\binom{(n-1)^2}{2}$ et la longueur du sous-sch\'ema support\'e par $Z$ est   $2n.\binom{n-1}{2}$, le sch\'ema restant est de longueur $\frac{n(n-1)(n-2)(n-3)}{2}$. Ce nombre est justement le degr\'e de la sous vari\'et\'e $\frak X_{n-1,n-2}\subset \P(H^{0}(O_{\P^{2}}(n-1))^{*})$
des courbes de degr\'e $(n-1)$ poss\`edant un point singulier d'ordre $(n-2)$ (ce que nous expliquerons plus en d\'etail  dans le paragraphe qui suivra l'exemple prototypique ci-dessous). Ceci prouve que la r\'eunion est disjointe.~$\Box$

\medskip

La proposition suivante explicite les liens entre les droites de saut de 
$E(Z)$ o\`u $Z$ est un groupe de points de longueur $2n$ en position lin\'eaire g\'en\'erale et celles de $E(Z\cup \{x\})$ o\`u $Z\cup \{x\}$ est  en position lin\'eaire g\'en\'erale.
%%%%%%%%%%%%%%
\begin{pro}
\label{fixe}
Soient $Z$ un groupe de points de longueur $2n$ en position lin\'eaire g\'en\'erale et 
$U_{Z}\subset \P^2$ l'ouvert des points $x$ tels que 
$Z\cup \lbrace x \rbrace$ soit en position lin\'eaire g\'en\'erale. 
On a alors,  
$$ S(E(Z))\subset \bigcap_{x\in U_{Z}} S(E(Z\cup \{x\})) $$ 
et surtout, lorsque $Z$ est  g\'en\'eral parmi les groupes de points du plan de longueur $2n$, leurs supports co\"incident. 
\end{pro}
%%%%%%%%%%%%%%%%%%%%%%%%%%%%%%
{\bf Remarque :} L'hypoth\`ese suppl\'ementaire sur $Z$ concernant l'\'egalit\'e des supports n'est pas n\'ec\'essaire \`a mon avis mais, sans elle, je n'ai pas de preuve valable.\\
{\bf Preuve de la proposition \ref{fixe}. }
Soit $x\in U_{Z}$, on consid\`ere  l'extension
(d'apr\`es $a)$ du th\'eor\`eme \ref{principal} par exemple ou bien \cite{DK} prop 2.9)
$$ \begin{CD}
0 @>>> E(Z) @>>>E(Z\cup \{x\})
@>>> O_{x^{\vee}}
@>>>0
\end{CD}
$$
Toute droite
de saut $l$ de $E(Z)$ est aussi une droite de saut de $E(Z\cup \{x\})$. En effet consid\'erons la restriction de le suite exacte ci-dessus \`a la droite $l$ on obtient
$$ \begin{CD}
0 @>>> O_{l}(n+1)\oplus O_{l}(n-2) @>>>E(Z\cup \{x\})\otimes O_{l}
@>>> O_{x^{\vee}\cap l}
@>>>0
\end{CD}
$$
L'injectivit\'e de la premi\`ere fl\`eche implique que $E(Z\cup \{x\})\otimes O_{l}\neq O_{l}(n)\oplus O_{l}(n)$ i.e que $l$ est une droite de saut pour $E(Z\cup \{x\})$, d'o\`u 
$$ S(E(Z))\subset \bigcap_{x\in U_{Z}} S(E(Z\cup \{x\})) $$
L'\'egalit\'e des supports sera prouv\'ee dans la proposition \ref{supports} du paragraphe qui suit.~$\Box$

\medskip

%%%%%%%%%%%%%%%%%%%%%%%%%%%%%%%%%%%%%%%%%%%%%%%%%%%%%%%%%%%%%%%
{\bf Exemple.} Soit $Z=\{ l_1,\cdots,l_8\}$. Le fibr\'e logarithmique
$E(Z)$ apr\`es
normalisation a pour classes de Chern $c_1=-1$ et $c_2=9$. Les droites de
saut $l_i$ sont des bisauteuses et
apparaissent avec une multiplicit\'e \'egale \`a $3$ (voir lemme \ref{singodd}), i.e la
longueur du sous-sch\'ema des droites de saut
support\'e par $Z$ est $24$. La longueur du sch\'ema des droites de saut
de $E(Z)$ est
$36$. Par ailleurs  ce fibr\'e \'etant d\'etermin\'e par la donn\'ee des
huit bisauteuses de $Z$ les douze
droites manquantes sont d\'etermin\'ees uniquement par la donn\'ee de $Z$.
Consid\'erons le pinceau de
cubiques passant par $Z$. Le lieu des points singuliers des cubiques du
pinceau sont les droites de saut de
$E(Z)$ distinctes de $Z$. Dans le $\P^9$ des cubiques du plan
l'hypersurface des cubiques
singuli\`eres, donn\'e par le discriminant, est de degr\'e $12$. Notre
pinceau de cubiques rencontre
cette hypersurface en douze points. 

\smallskip

Il est difficile de parler de cubiques passant par huit points sans dire un mot du neuvi\`eme.
%%%%%%%%%%%%%%%%%%%
\begin{rem}
Le pinceau de cubiques associ\'e aux huit points en position g\'en\'erale de $Z=\{
l_1,\cdots,l_8\}$ d\'efinit un neuvi\`eme point $l_9$. La droite correspondante n'est pas
une droite de saut de $E(Z)$.
\end{rem}
%%%%%%%%%%%
En effet, si une cubique du pinceau admet ce neuvi\`eme point comme point singulier
elle intersecte la cubique g\'en\'erale du pinceau en degr\'e $\ge 10$. Par cons\'equent
la cubique g\'en\'erale est d\'ecompos\'ee ce qui contredit l'hypoth\`ese de position
g\'en\'erale. 

\smallskip

Les cas limites sont les suivants : soit $C$ une conique lisse. Si $Z\in C$ alors 
$E(Z)$ est un fibr\'e de Schwarzenberger et son sch\'ema de droites de saut est le
diviseur 
$S(E(Z))=3C$ (voir \cite{Va1}, prop 2.1). Si $Z\setminus \{l_8\}\in C$, alors 
$S(E(Z))$ est support\'e par $C\cup \{l_8\}$.
%%%%%%%%%%%%%%%%%%
\section{Surfaces rationnelles, dualit\'e et droites de saut}
%%%%%%%%%%%%%%%%%%%%%%%%%%%%%%%%%%%%%%
A partir d'un groupe de points plan de longueur $2n$, on retrouve par des proc\'ed\'es classiques (\'eclatement, projections, vari\'et\'es tangentes...) les courbes mono\"idales associ\'ees (en ajoutant un point) ainsi que le groupe de points suppl\'ementaire. Avant de 
d\'etailler cette construction j'en donne ci-dessous la trame. 
\subsection{Squelette de la construction } 
Soit $Z_{0}$ un groupe de points plan en position g\'en\'erale de longueur $2n$. Les courbes de degr\'e $n$ avec un point singulier d'ordre $(n-1)$ passant par $Z_{0}$ forment, dans l'espace projectif des courbes de degr\'e $n$ passant par 
$Z_{0}$, une surface que nous noterons $S_{1}$.  Une section hyperplane de cette surface  par un hyperplan $(n-1)$-tangent en un point $x$  est une courbe birationnelle \`a la courbe mono\"{i}dale
$C(Z_0\cup \{x\})$.\\
Comme les courbes $C_{i}$ de degr\'e $(n-1)$  avec un point singulier $x_{i}$   d'ordre $(n-2)$ passant par $Z_{0}$ sont en nombre fini, \'egal \`a $N$, la surface $S_{1}$ contient  $N$ droites qui se contractent en un groupe de points $\Gamma=\lbrace x_{1}, \cdots, x_{N} \rbrace$ sur $\P^{2}$. En effet, \'etant donn\'e une courbe $C_{i}$ de degr\'e $(n-1)$ passant par $Z_{0}$ qui est $(n-2)$-singuli\`ere en un point $x_{i}$, le $\P^{1}$ de courbes de degr\'e $n$ obtenues en prenant la r\'eunion de $C_{i}$ et d'une droite passant par $x_{i}$ est bien contenu dans $S_{1}$.  Il en r\'esulte  que toutes les courbes mono\"{i}dales de la forme 
$C(Z_0\cup \{x\})$ ont $\Gamma$ pour lieu base. On obtient alors 
$S(E(Z_{0}))=Z_0^{n-1}\cup \Gamma$.
%%%%%%%%%%%%%%%%%%%%
\subsection{D\'etails de la construction}
%%%%%%%%%%%%%%%
Soit $V$ un espace vectoriel complexe de dimension $3$. On note $\P^2=\P V$ et $H^0(O_{\P^2}(1))=V$.  Consid\'erons
l'application canonique de d\'erivation  (application duale de la multiplication) $$\partial_{n-2} :
S^{n}V
\rightarrow S^{n-2}V \otimes S^{2}V $$
elle induit
l'homomorphisme suivant de fibr\'es
vectoriels
$$ \begin{CD}
0 @>>> S^{n-2}V^{*}\otimes O_{\P^{2}}(-2) @>\partial_{n-2}>>
S^{n}V^{*}\otimes O_{\P^{2}} @>>> E_{n-1}
@>>>0
\end{CD}
$$
o\`u la matrice $\partial_{n-2}$ est form\'ee des d\'eriv\'ees partielles
$(n-2)-$i\`emes d'une base de
$S^{n}V^{*}$. Le fibr\'e projectif $$\P(E_{n-1}) \subset \P^2\times
\P(S^{n}V^{*})$$  est d\'ecrit par l'incidence
$\{(x,f) ;  \partial_{n-2}f(x)=0\}$. La fibre au dessus d'un point
$x\in \P^2$ s'identifie \`a l'espace
projectif $\P(H^0({\frak m}_x^{n-1}(n)))$ o\`u ${\frak m}_x$ est le
faisceau d'id\'eaux  du point $x$.
Son image $\frak X_{n,n-1}$ dans
$\P(S^{n}V^{*})$ par la seconde projection est la vari\'et\'e des courbes
de degr\'e $n$ de $\P^2$ qui
s'annulent \`a l'ordre
$n-1$ en un point de $\P^2$, autrement dit les hyperplans $(n-1)$-tangents \`a la surface de Veronese $v_{n}(\P V)\subset \P S^{n}V$. Le fibr\'e projectif est 
birationnel \`a son image $\frak X_{n,n-1}$. La dimension, \'egale \`a  $2n+2$, de $\frak X_{n,n-1}$ ainsi que son degr\'e, \'egal \`a $\frac{(n+1)n(n-1)(n-2)}{2}$, se d\'eduisent des rang et classes de Chern du fibr\'e $E_{n,n-1}$. 

\smallskip

{\bf La surface $S_{1}$}
%%%%%%%%%%%%%%%%%%%%%%%%%%%%%%%%%%%%%%%

\smallskip

Soit $Z_0$ un groupe de points de $\P^2$ en position g\'en\'erale et de longueur $2n$. 
On consid\`ere la surface $S_{1}=\frak X_{n,n-1} \cap \P (H^0(J_{Z_0}(n))^{*})$. C'est une surface rationnelle provenant de l'\'eclatement de $\P^{2}$ le long d'un sch\'ema de longueur finie
$Z_1$ d\'efini par $Z_{0}$ que nous explicitons ci-dessous.

\medskip

La restriction
$ \begin{CD}
 (S^{n}V)^{*} @>res>> H^0(J_{Z_0}(n))^{*}
\end{CD}
$ compos\'ee avec l'homorphisme
 $\partial_{n-2}$ induit une suite exacte
$$ \begin{CD}
0 @>>> S^{n-2}V^{*}\otimes O_{\P^{2}}(-2) @>M>>H^0(J_{Z_0}(n))^{*}\otimes
O_{\P^{2}} @>>> J_{Z_1}(n(n-1))
@>>>0
\end{CD}
$$
o\`u le sch\'ema de longueur finie $Z_{1}$   est le lieu base du syst\`eme lin\'eaire d\'efini par les mineurs maximaux
de la matrice $M$. Le support de $Z_1$ est form\'e des points $x_{1}$ du plan pour lesquels $ h^0(J_{Z_0}\otimes {\frak
m}_{x_1}^{n-1}(n))\ge 2$. C'est le cas par exemple lorsque $x_1\in Z_0$, ce qui
montre  que $Z_0\subset Z_1$. La surface $S_{1}$ est l'image de l'\'eclat\'e de $\P^{2}$ le long de $Z_{1}$ par le morphisme compos\'e :
$$  \P(J_{Z_1}(n(n-1)))\hookrightarrow \P V\times \P(H^0(J_{Z_0}(n))^{*})
\stackrel{pr_2}\rightarrow
\P(H^0(J_{Z_0}(n))^{*})$$

\smallskip

{\bf Les courbes mono\"idales}
%%%%%%%%%%%%%%%%%%%%%%%%%%%%%%%%%%%%%%%

\smallskip

Etant donn\'e $x\in \P^{2}\setminus Z_{0}$, le conoyau de l'application compos\'ee $$S^{n-2}V^{*}\otimes O_{\P^{2}}(-2)\rightarrow H^0(J_{Z_0}(n))^{*}\otimes
O_{\P^{2}} \rightarrow 
H^0(J_{Z_0\cup \{x\}}(n))^{*}\otimes O_{\P^{2}} $$ est  support\'e par la courbe mono\"{i}dale
$C(Z_0\cup \{x\})$. Cette courbe est une section de $J_{Z_1}(n(n-1))$ birationnelle \`a la courbe
d'intersection $\frak X_{n,n-1} \cap \P (H^0(J_{Z_0\cup \{x\}}(n))^{*})$, elle est donn\'ee
ensemblistement par $\{y\in \P^2 \mid h^0(J_{Z_0\cup \{x\}}\otimes \frak m_y^{n-1}(n))\neq 0\}$.

\smallskip

{\bf Le groupe de point \'eclat\'e}
%%%%%%%%%%%%%%%%%%%%%%%%%%%%%%%%%%%%%%%

\smallskip

On voudrait maintenant donner une description plus pr\'ecise du groupe de
points $Z_1$. Compte tenu de ce que nous avons vu ci-dessus, aux points de $Z_1$ il
appara\^{\i}t au moins un $\P^1$ de courbes de degr\'e $n$ qui sont $(n-1)$-singuli\`eres. On
veut montrer que cette situation (pour $Z_{0}$ en position g\'en\'erale) est enti\`erement d\'ecrite par  l'union d'une courbe de degr\'e $(n-1)$ ayant un point singulier d'ordre $(n-2)$
et d'une droite passant par ce point singulier. 
%%%%%%%%%%%
\begin{lem}
\label{pinceausing}
Soit $x\in \P^2\setminus Z_0$ un point tel que $h^0(J_{Z_0}\otimes {\frak
m}_x^{n-1}(n))\ge 2$. Alors une section g\'en\'erale  de $H^0(J_{Z_0}\otimes {\frak
m}_x^{n-1}(n))$ est la r\'eunion d'une courbe de degr\'e $(n-1)$ qui est $(n-2)$ singuli\`ere en
$x$ et d'une droite passant par $x$.
\end{lem}
%%%%%%%%%
{\bf Preuve du lemme \ref{pinceausing}.}
%%%%%%%%%%%
Un pinceau de courbes
singuli\`eres poss\`ede une composante commune. Consid\'erons deux courbes du pinceau. Si elles n'ont pas de composante commune elles se coupent le long d'un sch\'ema de longueur 
$n^2$ qui contient un sous-sch\'ema de longueur au moins $(n-1)^2$ au point base singulier. Ceci implique qu'en dehors du point singulier elles doivent se couper au plus le long de $2n-1$ points, ce qui contredit l'hypoth\`ese.  Montrons que cette composante
est une courbe de ${\frak X_{n-1,n-2}}$ singuli\`ere en $x$.  

\smallskip

Comme $Z_{0}$ est en position g\'en\'erale le lieu
$\{x\in \P^{2}, h^0({\cal J}_{Z_{0}}\otimes {\frak m}_x^{n-3}(n-2))\neq 0 \}$ est vide et le lieu $\{x\in \P^{2}, h^0({\cal J}_{Z_{0}}\otimes {\frak m}_x^{n-2}(n-1))\neq 0 \}$
est un sch\'ema de longueur finie. Par cons\'equent une courbe de 
$S_{1}$ (c'est \`a dire une courbe de $\P^{2}$ correspondante \`a un point de $S_{1}$) qui n'est pas irr\'eductible est obligatoirement de la forme $C\times L$ avec 
$C\in {\frak X_{n-1,n-2}}\cap \P(H^0(J_{Z_0}(n-1))^{*})$ et $L$ une droite passant par $x$.~$\Box$

\medskip

Comme ${\frak X_{n-1,n-2}}\cap \P(H^0(J_{Z_0}(n-1))^{*})$ est
un  groupe de points de longeur $N={\rm deg}{\frak X_{n-1,n-2}}$ (les courbes 
$C_{i}$ qui sont  $(n-2)$-singuli\`eres aux points
$x_{i}$)
la surface
$S_1={\frak X_{n,n-1}}\cap \P(H^0(J_{Z_0}(n))^{*})$ contient alors  
${\rm deg}(
\frak X_{n-1,n-2})$ droites. Ces
droites sont contract\'ees sur $\P^2$ (le m\^eme point singulier pour toute
la droite). On note $\Gamma=\lbrace x_{1}, \cdots, x_{N}\rbrace $
le groupe de points image dans $\P^2$. 
~$\Box$
%%%%%%%%%%%
\begin{pro}
\label{supports} 
Le support de $\bigcap_{x\in U_{Z_{0}}}C(Z_{0}\cup \lbrace x \rbrace)$, $Z_{1}$ et de $S(E(Z_{0}))$ est $Z_0\cup \Gamma$. Plus pr\'ecis\'ement, 
$Z_1=\tilde{Z_0}\cup \Gamma$  et 
$S(E(Z_{0}))=Z_{0}^{n-1}\cup \Gamma$
o\`u $\tilde{Z_0}$ est support\'e
par $Z_0$ et est localement
intersection compl\`ete $(n-1)\times (n-1)$.
\end{pro}
%%%%%%%%%%%
{\bf Preuve de la proposition \ref{supports}.} D'apr\`es le lemme 
 pr\'ec\'edent  $Z_1$ est  ensemblistement la
r\'eunion de $Z_0$ et de $\Gamma$. Par
ailleurs 
lorsque $Z_{0}$ est en position g\'en\'erale les courbes mono\"idales 
$C(Z_{0}\cup \lbrace x \rbrace)$ engendrent l'espace $\P(H^{0}J_{Z_{1}}(n(n-1)))$. Autrement dit $Z_{1}=\bigcap_{x\in U_{Z_{0}}}C(Z_{0}\cup \lbrace x \rbrace)$. On en d\'eduit que 
$S(E(Z_{0}))\subset Z_{1}$ ou encore que $Z_{1}$ contient le 
$(n-2)$-i\`eme voisinage
infinit\'esimal de $Z_0$. Comme $Z_1$ est localement intersection
compl\`ete (c'est le lieu d'annulation d'une section d'un fibr\'e vectoriel de rang deux!) et
comme
$l(O_{Z_1})= (n-1)^2l(O_{Z_0})+ {\rm deg}( \frak X_{n-1,n-2})$  on en d\'eduit le
r\'esultat.~$\Box$
%%%%%%%%%%%%%%%%%%%%%%%%%%%%%%

\end{document}